\documentclass[psamsfonts]{amsart}
\usepackage[utf8]{inputenc}
\usepackage[english]{babel}
\usepackage[T1]{fontenc}
\usepackage{amsfonts}
\usepackage{amsmath}
\usepackage{amsfonts}
\usepackage{bbm}
\usepackage{graphicx}
\usepackage{array}
\usepackage{lipsum}
\usepackage{amsaddr}
\usepackage{subcaption}

\usepackage{adjustbox}
\usepackage{transparent}

\usepackage{isomath}
\usepackage{mathtools}
\usepackage{mathrsfs}
\usepackage{amsthm}
\usepackage{slashed}
\usepackage{amssymb}
\usepackage{enumitem}
\usepackage{hyperref}
\hypersetup{colorlinks=true, 
citecolor=red}
\usepackage[new]{old-arrows}

\theoremstyle{definition}
\newtheorem{mdef}{{Definition}}[section]

\theoremstyle{definition}
\newtheorem{mex}{Example}[section]

\theoremstyle{definition}
\newtheorem{mrmk}{{Remark}}[section]

\theoremstyle{plain}
\newtheorem{mth}{Theorem}[section]

\theoremstyle{plain}
\newtheorem{mlem}{{Lemma}}[section]

\theoremstyle{plain}
\newtheorem{mprop}{{Proposition}}[section]

\theoremstyle{plain}
\newtheorem{mcor}{{Corollary}}[section]

\theoremstyle{definition}

\newlength{\bibitemsep}
\newlength{\bibparskip}
\let\oldthebibliography\thebibliography
\renewcommand\thebibliography[1]{%
  \oldthebibliography{#1}%
  \setlength{\parskip}{\bibitemsep}%
  \setlength{\itemsep}{\bibparskip}%
}

\usepackage{tikz}
\usepackage{pgfplots}
\usetikzlibrary{arrows.meta}
  \pgfplotsset{compat=newest}
  \usetikzlibrary{plotmarks}
  \usepackage{grffile}

\pgfplotsset{every axis/.append style={
        scaled ticks = false, 
        tick label style={/pgf/number format/fixed}
    }
}

\usetikzlibrary{external}
\newlength\figureheight 
\newlength\figurewidth 

\usepackage{shellesc}
\usetikzlibrary{plotmarks}
\usetikzlibrary{patterns}
\usepackage{tikz-3dplot}
\usepackage{multirow}
\usepackage{caption}
\usepackage{subcaption}

\DeclareMathSymbol{\upLambda}{\mathalpha}{operators}{3}

\title[Noncommutative Laplacian]{Noncommutative Laplacian and numerical approximation of Laplace-Beltrami spectrum of compact Riemann surfaces}
\author{Damien Tageddine, Jean-Christophe Nave}
\address{Department of Mathematics and Statistics, McGill University}
\begin{document}
\maketitle
\begin{abstract}
\noindent
We derive a numerical approximation of the Laplace-Beltrami operator on compact Riemann surfaces embedded in $\mathbb{R}^3$ with an effective $S^1$ action. To do so, we use a noncommutative Laplace operator defined on the space of finite dimensional hermitian matrices. This operator is derived from a foliation of the surface obtained under the $S^1$-action. We present numerical results in the case of an ellipsoid and surfaces of genus $g\geq 1$.
\end{abstract}
{
  \hypersetup{linkcolor=black}
  \tableofcontents
}
\section{Introduction}
In the present work, we are interested in the eigenvalue problem to determine the spectrum of $-\mathrm{\Delta}_g$. The determination of the spectrum of the Laplace operator \cite{canzani_analysis_nodate, rosenberg_laplacian_1997, jakobson_how_2005} is central in several areas of mathematics such as spectral analysis, mathematical physics and geometry processing \cite{Reuter2006LaplaceBeltramiSA}. Since closed-form solutions are rarely available except for very simple geometry such as the sphere or the flat torus, numerical methods are typically used. These methods vary based on the representation of the surface, e.g. triangulated mesh, point cloud or implicit surfaces.
\subsection{Main results}
We introduce a numerical approximation of the Laplace Beltrami operator on compact Riemann surfaces based on a noncommutative analogue of the Laplacian. The construction relies on replacing the commutative algebra of smooth functions with a sequence of finite-dimensional matrix algebras, in the spirit of matrix regularization \cite{SHIMADA2004297, hoppe_quantization_2021, arnlind_graph_2008, arnlind_multi_2010, arnlind_construction_2022}. Within this framework, we define a discrete Laplace operator built from commutators of quantized coordinate functions. We prove the following convergence theorem.
\begin{mth}
Let $(\mathrm{\Sigma},\omega)$ be a compact toric surface embedded in $\mathbb{R}^3$ equipped with the induced metric tensor $g$. There exists a sequence of maps $(T_\hbar,\hbar)$ such that the eigenmatrix of $\mathrm{\Delta}_\hbar$ converges to an eigenfunction of $\mathrm{\Delta}$ with eigenvalue given by the limit $\lambda=\lim \lambda_\hbar$ as $\hbar\rightarrow 0$.
\end{mth}
\noindent
Our main result shows that the eigenmatrices of this noncommutative Laplacian converge to eigenfunctions of the classical Laplace–Beltrami operator, with convergence of the corresponding eigenvalues. This provides a structure-preserving discretization scheme that is particularly well-suited for surfaces of arbitrary genus, admitting $S^1$-actions, but which can also be extended to more general settings. We further demonstrate the applicability of this method by numerically computing eigenvalue spectra of the Laplacian on the sphere, an ellipsoid, the standard torus and a double-torus. We compare our results with existing literature \cite{volkmer_laplace-beltrami_2021},\cite{eswarathasan_laplace-beltrami_2021},\cite{volkmer_eigenvalues_2023}, and \cite{volkmer_laplace-beltrami_2024} respectively.
\subsection{Preliminaries}
Let $\mathrm{\Sigma}$ be a surface with local coordinates $u^1$ and $u^2$, one can define 
\begin{equation}
\frac{1}{\sqrt{|g|}}\left\lbrace f,g \right\rbrace  =\frac{1}{\sqrt{|g|}}\left( \frac{\partial f}{\partial u^1}\frac{\partial h}{\partial u^2}-\frac{\partial h}{\partial u^1}\frac{\partial f}{\partial u^2}\right) 
\end{equation}
where $g$ is the determinant of the induced metric tensor, and one readily checks that $(C^\infty(\mathrm{\Sigma}),\left\lbrace \cdot , \cdot \right\rbrace )$ is a Poisson algebra. The Laplace-Beltrami operator $\mathrm{\Delta}_g$ acts on smooth functions by
\begin{equation}
\mathrm{\Delta}_gf=\frac{1}{\sqrt{|g|}}\partial_i\left( \sqrt{g}g^{ij}\partial_jf\right) 
\end{equation}
The eigenvalue problem given by
\begin{equation}
\mathrm{\Delta}_gf=\lambda f
\end{equation}
defines the spectrum of the operator $\mathrm{spec}(\mathrm{\Delta}_g)$ on $(\mathrm{\Sigma},g)$. Since this spectrum is deeply tied to the geometry of the surface, approximating it is a central problem in spectral geometry and numerical analysis.\\

\noindent
Following the noncommutative approach \cite{connes_noncommutative_1994, klimek_quantum_1992, klimek_quantum_1992-1}, instead of discretizing the manifold directly, we encode the geometry into a noncommutative Poisson algebra and approximate smooth functions by Hermitian matrices \cite{tageddine_statistical_2023, tageddine_structure_2024, tageddine_noncommutative_2025}.
This procedure is referred to as \textit{quantization} or, in the present context, \textit{matrix regularization}, whereby the infinite-dimensional algebra of smooth functions on the surface is approximated by an algebra of $N\times N$ matrices. We begin by providing an informal description of the matrix quantization method and refer the reader to Section \ref{Section3}, Definition \ref{matrixdef} for the formal definition.
For each positive integer $N$, we construct a linear map 
\begin{equation*}
    T_N:C^\infty(\mathrm{\Sigma},\mathbb{R}) \rightarrow \mathrm{Herm}(N)
\end{equation*}
where $\mathrm{Herm}(N)$ denotes the space of Hermitian $N\times N$  matrices. Let $\hbar(N)$ be a real-valued strictly decreasing function such that $\lim_{N\rightarrow \infty}N\hbar(N)<\infty$. The principal properties required of the regularization map $T_N$ are
\begin{align*}
\lim_{N\rightarrow\infty} \|T_N(f)T_N(g)-T_N(fg)\|=0,\quad \lim_{N\rightarrow \infty}\left| \left|\frac{1}{i\hbar_N}\left[ T_N(f),T_N(g)\right] - T_N(\left\lbrace f,g \right\rbrace )\right| \right|=0
\end{align*}
If such family of maps $\left\lbrace T_N\right\rbrace_{N\in\mathbb{N}}$ exists and satisfies the conditions above, then one speaks of a matrix regularization of the Poisson algebra $(C^\infty(\mathrm{\Sigma},\mathbb{R}),\left\lbrace \cdot,
\cdot\right\rbrace)$.
\subsection{Related results}
There are a number of relevant approaches to approximate $\mathrm{spec}(-\mathrm{\Delta}_g)$ over a curved surface that require discussion. We begin by mentioning the work \cite{dziuk_finite_2013} on finite element methods for surface PDEs. Here, the authors use of surface finite elements to compute solutions to the Poisson problem for the Laplace-Beltrami operator on curved surfaces. This approach relies on either triangulated surfaces or implicit degenerate equations, on which weak forms of the given surface elliptic operator are defined.\\
Another method of relevance is the closest point method \cite{macdonald_level_2008,macdonald_implicit_2010}. This numerical technique for solving PDEs on surfaces consist in embedding the given surface into a higher-dimensional Euclidean space. The key idea is to extend functions from the surface to the surrounding space by assigning to each point in the embedding space the value of the function at its closest point on the surface. This allows one to use standard Cartesian grid discretizations while naturally enforcing the surface geometry. In particular, the method can be used to compute eigenvalues of the Laplace–Beltrami operator by discretizing the extended problem in the embedding space and applying standard eigensolvers to approximate the intrinsic spectral properties of the surface. The closest point method can be categorized into the broader framework of point clouds methods: techniques such as the meshless Laplace Beltrami or diffusion maps approximate the operator from sampled points without explicit triangulation \cite{belkin_constructing_2009, belkin_towards_2008}.\\
Closer to our approach, we mention the works on quantization and matrix regularizations applied to discretization of compact Riemann surfaces \cite{sykora_fuzzy_2017,arnlind_fuzzy_2009,adachi_matrix_2020,arnlind_low_2020,arnlind_noncommutative_2007}. For instance, the Berezin-Toeplitz quantization provides finite dimensional operator models of function algebras. Furthermore, these methods have been developed to derive structure-preserving  \cite{tageddine_structure_2024}. Indeed, recent advances have also focused on the relation between noncommutative geometries and discrete differential structures \cite{tageddine_noncommutative_2025, tageddine_statistical_2023}. This is in particular exemplified by the Berezin-Toeplitz quantization. This approach discretizes certain types of Poisson algebras while maintaining the Lie-algebra structure in the transition from continuous to discrete systems \cite{modin_eulerian_2024,modin_liepoisson_2020}. Hence, quantization provides a robust framework for discrete geometry applications.  We also mention fundamental advances on geometric reconstruction from finite dimensional spectral triples and point cloud dynamics \cite{nielsen_dynamical_2023}.
\subsection{Outline of the paper}
The paper is organized as follows. In Section \ref{Section2}, we introduce the geometric framework underlying our construction, namely the foliation of compact Riemann surfaces induced by an $S^1$-action, and the associated momentum map. Section \ref{Section3} presents the matrix regularization of the Laplace–Beltrami operator on a manifold defined with a single coordinate chart. We first recall the necessary preliminaries on Poisson brackets and the Laplace–Beltrami operator, and then construct the noncommutative analogue through a quantization map defined on matrix algebras. Section \ref{Section4}  treats the general case of a manifold given by several coordinate charts. Section \ref{Section5} is devoted to numerical experiments that illustrate the convergence and robustness of the proposed method on the sphere, the ellipsoid and the embedded torus $\mathbb{T}^2\subset \mathbb{R}^3$.
\section{Momentum map, Foliation and $S^1$-action}
\label{Section2}
Consider a $2$-dimensional compact connected symplectic manifold $(M,\omega)$. In addition, let us assume that $S^1$ acts smoothly on $M$ i.e. we have a group homomorphism:
\begin{equation}
\alpha:\begin{array}{c}
S^1\rightarrow \mathrm{Diff}(M) \\ 
g\mapsto \alpha_g
\end{array} 
\end{equation}
such that the evaluation map $\mathrm{ev}_\alpha(p,g)=\alpha_g(p)$ is smooth at any point $p\in M$. We will also suppose that the action of $S^1$ is symplectic:
\begin{equation}
\alpha:S^1\rightarrow \mathrm{Sympl}(M,\omega)\subset \mathrm{Diff}(M)
\end{equation}
is a symplectomorphism and hamiltonian if there exists a function $H:M\rightarrow \mathbb{R}$ with $dH=\iota_X\omega$, where $X$ is the vector field generated by $\alpha$. Such function is called a \textit{momentum map}.\\
\noindent
A particular class of symplectic manifolds satisfying the conditions above are toric manifold. First, let us recall that an action of a group $G$ on a manifold $M$ is called \textit{effective} if each group element $g\neq e$ moves at least one $p\in M$, that is,
\begin{equation*}
\cap_{p\in M}G_p=\left\lbrace e \right\rbrace 
\end{equation*}
where $G_p$ is the stabilizer of $p$.
\begin{mdef}[Toric manifold \cite{da_silva_lectures_2008}]
A (symplectic) toric manifold is a compact connected symplectic manifold $(M,\omega)$ equipped with an effective hamiltonian action of a torus $\mathbb{T}$ of dimension
\begin{equation}
\dim \mathbb{T}=\frac{1}{2}\dim M
\end{equation}
and a choice of a corresponding moment map.
\end{mdef}
\noindent
Let us focus on an effective action of $S^1$ on a symplectic manifold. In that case, a toric manifold is a surface $M$ with an effective action of $S^1$ with its associated hamiltonian map $H$.
\begin{mex}
The unit sphere $\mathbb{S}^2$ with it standard $2$-form $\omega=d\theta\wedge dz$ (cylindrical coordinates) with vector field $X=\frac{\partial}{\partial\theta}$. Each orbital is a horizontal circle $\left\lbrace (\theta + t,z):t\in \mathbb{R}\right\rbrace $ which is given by an action of $S^1$: 
\begin{equation}
\begin{array}{c}
\alpha: S^1\rightarrow \mathrm{Sympl}(\mathbb{S}^2,\omega) \\ 
t\mapsto \text{rotation by angle $t$ around $z$-axis}
\end{array} 
\end{equation}
This action is hamiltonian with moment map given by the height of a point:
\begin{equation}
h(\theta,z)=z
\end{equation}
and also effective making $(\mathbb{S}^2,\omega)$ into a toric manifold.
\end{mex}
\noindent
In the general case of a toric surface $\mathrm{\Sigma}$ embedded in $\mathbb{R}^3$ with an action of $S^1$, the height function
\begin{equation}
h:\mathbb{R}^3\rightarrow \mathbb{R} \quad \mu(r,\theta,z)=z
\end{equation}
is the associate moment map. In addition, it defines a Morse function; that is, a smooth real-valued function on $\mathrm{\Sigma}$ with no degenerate critical points (i.e. a critical point at which the Hessian matrix is singular). Moreover, from this Morse function $h$, we obtain a foliation of $\mathrm{\Sigma}$; that is a decomposition into submanifolds. Each of these submanifolds is called a \textit{leaf} of the foliation, and one says that $\mathrm{\Sigma}$ is \textit{foliated by the leaves}. The leaves are given by the inverse image $h^{-1}(\lambda)$ for $\lambda\in \mathbb{R}$. If $\lambda$ is a regular value of $h$, then the leaves are one the union of one dimensional closed curves (homeomorphic to $S^1$); one speaks of a \textit{regular leaf}. If $\lambda$ is an (isolated) singular value, then the dimension of leaves is allowed to vary and is either a point or a closed curve. One speaks instead of a \textit{singular leaf}.\\
A decomposition of $\mathrm{\Sigma}$ is then given as a disjoint union of connected subsets $V_\alpha$ of leaves,
\begin{equation}
M=\bigcup_{\lambda}V_\lambda
\end{equation}
where each leaf $V_\lambda$ may be defined as a level set of the height function $h$:
\begin{equation}
V_\lambda=h^{-1}(\lambda)\qquad \lambda\in \mathbb{R}.
\end{equation}
\noindent
In addition, every regular level $h^{-1}(\lambda)$ is a $1$-dimensional compact submanifold of $\mathrm{\Sigma}$, hence:
\begin{equation}
h^{-1}(\lambda) \simeq \sqcup_n S^1.
\end{equation}
\begin{figure}[htp!]
\centering
    \begin{subfigure}[t]{0.5\textwidth}
        \centering
		\setlength\figureheight{0.200\textheight} 
		\setlength\figurewidth{0.8\linewidth} 
        \includegraphics[scale=0.4]{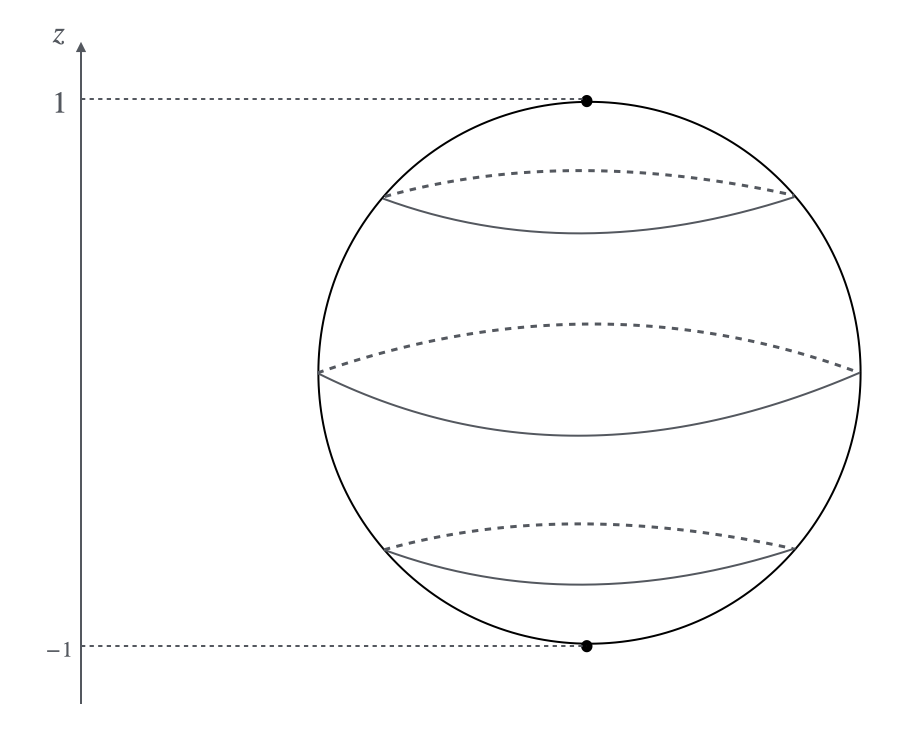}
        \caption{the unit sphere}
    \end{subfigure}
    \begin{subfigure}[t]{0.4\textwidth}
        \centering
		\setlength\figureheight{0.200\textheight} 
		\setlength\figurewidth{0.8\linewidth} 
        \includegraphics[scale=0.4]{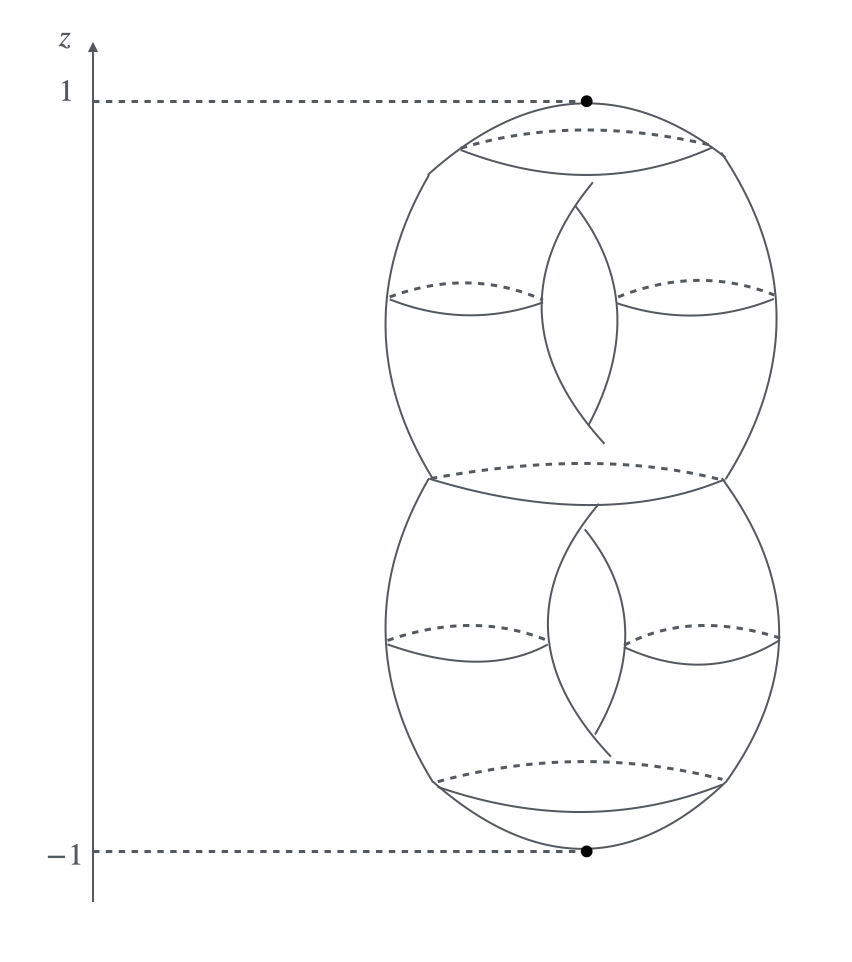} 
        \caption{the double torus}
    \end{subfigure}
    \caption{Foliation of embedded compact surfaces in $\mathbb{R}^3$}
\end{figure}
Consider a function $f:\mathbb{R}^3\rightarrow \mathbb{R}$ restricted on a leaf $\phi^{-1}(\lambda)$ can be written as the direct sum
\begin{equation}
f=(f_1,f_2,\cdots, f_n) \quad \text{s.t.} \quad f_i:[-1,1]\times S^1 \rightarrow \mathbb{R}.
\end{equation}
We can then consider functions $f$ such that for each $z\in [-1;1]$, $f_i(z,\cdot)\in L^2(S^1)$ and each $\theta \in S^1$, $f_i(\cdot,\theta)\in C([-1,1])$. We can write $f_i$ using its Fourier expansion
\begin{equation}
f_i(z,\theta)=\sum_{n\in \mathbb{Z}}a_n(z)e^{in\theta}.
\end{equation}
The goal of the next section is to construct a quantization map, which is a linear map such that $T_\alpha(f)$ is a hermitian map. This will be done using truncated Fourier series of the function $f$.
\section{Matrix regularization of the Laplace-Beltrami operator}
\label{Section3}
\subsection{Preliminaries}
Consider a surface $\mathrm{\Sigma}$ embedded in $\mathbb{R}^3$ defined by the embedding coordinates 
\begin{equation}
\mathrm{\Sigma}\rightarrow \mathbb{R}^3 \qquad (\theta,z)\mapsto (x(z,\theta),y(z,\theta),z)
\end{equation}
This manifold has the natural Poisson-bracket
\begin{equation}
\frac{1}{\sqrt{|g|}}\left\lbrace f,h \right\rbrace = \frac{1}{\sqrt{|g|}}\left( \partial_{\varphi}f\partial_zh-\partial_zf\partial_\varphi h\right) 
\end{equation}
where $|g|$ is the determinant of the induced metric tensor, and one can readily see that $(C^\infty(\mathrm{\Sigma}),\left\lbrace \cdot,\cdot \right\rbrace)$ is a Poisson algebra. The pullback metric on $\mathrm{\Sigma}$ can be expressed using the embedding coordinates:
\begin{equation}
g_{ab}=\sum_{i=1}^3\frac{\partial x^i}{\partial \sigma_a}\frac{\partial x^i}{\partial \sigma_b}
\end{equation}
and one can show that the square-root of the determinant of the metric tensor can be obtained as 
\begin{equation}
\sqrt{|g|}=\sqrt{\sum_{i>j}\left\lbrace x_i,x_j \right\rbrace^2 }.
\label{determetric}
\end{equation}
Similarly, the Laplace-Beltrami operator on $\mathrm{\Sigma}$ can also be expressed in terms of the Poisson bracket.
\begin{mlem}
On the Poisson manifold $(\mathrm{\Sigma},\left\lbrace\cdot,\cdot \right\rbrace )$, the Laplace-Beltrami operator can be written as
\begin{equation*}
\mathrm{\Delta}_g(f)=\sum_{i=1}^3\frac{1}{\sqrt{|g|}}\left\lbrace x^i,\frac{1}{\sqrt{|g|}}\left\lbrace x^i, f \right\rbrace \right\rbrace 
\end{equation*}
for any $f\in C^\infty(\mathrm{\Sigma})$.
\end{mlem}
\begin{proof}
Starting with the expression of the Laplace-Beltrami operator in local coordinates 
\begin{equation*}
\mathrm{\Delta}_g(f)=\frac{1}{\sqrt{|g|}}\partial_r\left( \frac{1}{\sqrt{|g|}}gg^{rs}\partial_sf\right) 
\end{equation*}
we use the fact that $gg^{rs}=\varepsilon^{ru}\varepsilon^{sv}g_{uv}$, then we can write
\begin{equation*}
\frac{1}{\sqrt{|g|}}\partial_r\left( \frac{1}{\sqrt{|g|}}gg^{rs}\partial_sf\right) = \frac{1}{\sqrt{g}}\partial_r\left( \frac{1}{\sqrt{g}} \varepsilon^{ru}\varepsilon^{sv}g_{uv}\partial_sf\right) 
\end{equation*}
Moreover, let us recall that $g_{uv}=\partial_ux^j\partial_vx^j$ which leads to
\begin{align*}
\mathrm{\Delta}_g(f)= \frac{1}{\sqrt{|g|}}\partial_r\left( \frac{1}{\sqrt{|g|}} \varepsilon^{ru}\varepsilon^{sv}\partial_ux^j\partial_vx^j\partial_sf\right).
\end{align*} 
Using the Poisson bracket to collect terms in the right-hand side, we get:
\begin{equation*}
\frac{1}{\sqrt{|g|}}\partial_r\left( \frac{1}{\sqrt{|g|}}gg^{rs}\partial_sf\right)=\frac{1}{\sqrt{|g|}}\partial_r\left(\frac{1}{\sqrt{|g|}}\varepsilon^{ru}\left\lbrace x^i,f \right\rbrace \partial_u x^j  \right).
\end{equation*}
Expanding the right-hand side and again collecting terms in the Poisson bracket, we obtain
\begin{equation*}
\mathrm{\Delta}_g(f)=\frac{1}{\sqrt{|g|}}\left\lbrace x^i, \frac{1}{\sqrt{|g|}}\left\lbrace x^i,f \right\rbrace \right\rbrace + \frac{1}{|g|}\left\lbrace x^i,f \right\rbrace \varepsilon^{ru}\partial^2_{ru}x^j.
\end{equation*}
The second term vanishes since $\varepsilon^{ru}$ is antisymmetric and $\partial^2_{ru}x^j$ is symmetric in the two indices.
\end{proof}
\subsection{Matrix quantization}
\begin{mdef}[Matrix quantization]
\label{matrixdef}
Let $N_1,N_2,\dots$ be a strictly increasing sequence of positive integers, let $\left\lbrace T_\alpha\right\rbrace $ for $\alpha=1,2,\dots$ be linear maps from $C^\infty(\mathrm{\Sigma},\mathbb{R})$ to hermitian $N_\alpha\times N_\alpha$ matrices and let $\hbar(N)$ be a real-valued strictly decreasing function such that $\lim_{N\rightarrow \infty}N\hbar(N)<\infty$.\\
Consider a linear map $T_\alpha$ from smooth functions to hermitian $N_\alpha\times N_\alpha$ matrices, the main properties of the regularization are
\begin{align}
&\lim_{\alpha\rightarrow\infty} \|T_\alpha(f)\|<\infty, \label{ax1} \\
&\lim_{\alpha\rightarrow\infty} \|T_\alpha(f)T_\alpha(g)-T_\alpha(fg)\|=0, \label{ax2}\\
&\lim_{\alpha\rightarrow \infty}\left| \left|\frac{1}{i\hbar_\alpha}\left[ T_\alpha(f),T_\alpha(g)\right] - T_\alpha(\left\lbrace f,g \right\rbrace )\right| \right|=0. \label{ax3}
\end{align}
\end{mdef}
\noindent
Let $(T_\alpha,\hbar_\alpha)$ be a matrix regularization of $(\mathrm{\Sigma},\omega)$, where $\mathrm{\Sigma}$ is embedded in $\mathbb{R}^m$ via coordinates $(x_1,\dots,x^m)$ with a symplectic form 
\begin{equation}
\omega=\rho(u^1,u^2)du^1\wedge du^2
\end{equation}
inducing the Poisson bracket $\left\lbrace f,g\right\rbrace =\frac{1}{\rho}\varepsilon^{ab}\partial_af\partial_bg$. Furthermore, we let $\left\lbrace \gamma_\alpha\right\rbrace $ be a convergent sequence converging to $\gamma=\sqrt{g}/\rho$ and we set 
\begin{equation}
T_\alpha(x^i)=X^i_\alpha.
\end{equation}
The \textit{noncommutative Laplacian} on $\mathrm{\Sigma}$ is a sequence $\left\lbrace \mathrm{\Delta}_\alpha\right\rbrace $ of linear maps defined as 
\begin{equation}
\mathrm{\Delta}_\alpha = -\frac{1}{\hbar_\alpha^2}\gamma^{-1}_\alpha [X_i,\gamma^{-1}_\alpha[X_i,\cdot ] ]
\end{equation}
 where $\gamma_\alpha$ can be taken to be 
 \begin{equation}
 \gamma_\alpha = \sqrt{\sum_{i>j}-\left( \frac{[X^i_\alpha,X^j_\alpha]}{\hbar_\alpha}\right)^2 }.
 \end{equation}
 Consider now a surface $\mathrm{\Sigma}$ embedded in $\mathbb{R}^3$ defined by the embedding coordinates 
\begin{equation}
\mathrm{\Sigma}\rightarrow \mathbb{R}^3 \qquad (\theta,z)\mapsto (x(z,\theta),y(z,\theta),z)
\end{equation}
Such a surface possesses an $S^1$-action leaving the $z$-coordinate invariant and for which a momentum map is given by
\begin{equation}
\mu:\mathbb{R}^3\rightarrow \mathbb{R} \qquad \mu(x,y,z)=z.
\end{equation}
As explained in Section \ref{Section2}, from this action one can get a foliation of $\mathrm{\Sigma}$ using $\mu$ as a Morse function following the previous section. Hence, the restriction of functions on $\mathrm{\Sigma}$ to a leaf $\mu^{-1}(\lambda)$ admits a Fourier decomposition. In particular, we can write the formal Fourier series of the coordinate functions:
\begin{align*}
x(z,\theta)=\sum_{n\in \mathbb{Z}}x_n(z)e^{in\theta}\qquad y(z,\theta)=\sum_{n\in \mathbb{Z}}y_n(z)e^{in\theta}
\end{align*}
along with their truncations. We will now assume some regularity on the Fourier expansion and consider that all the functions considered from now on have a restriction in $C^\infty(S^1)\cap L^2(S^1)$. In particular, for any such function, we have:
\begin{equation}
-i\frac{\partial^\alpha f}{\partial\theta^\alpha}=\sum_{n\in\mathbb{Z}} n^\alpha f_n(z)e^{in\varphi}.
\label{fouriercoeff}
\end{equation}
In order to define a quantization map, we also discretize the $z$ coordinate for $z\in [a,b]$: for $ N>0$,
\begin{equation}
z(n)=a+\frac{(b-a)\beta}{N}n, \qquad 1\leq n\leq N-1
\end{equation}
where $\beta>0$ is some parameter to be chosen. We also define the two-variable function:
\begin{equation}
z(n,m)=z\left( \frac{n+m}{2}\right) = a+\frac{(b-a)\beta}{2N}(n+m), \qquad 1\leq n,m\leq N-1.
\end{equation}
Notice that in particular,
\begin{equation}
z(n,n)=z(n).
\end{equation}
In order to simply the notations, we denote by $\hbar_N$ (or simply $\hbar$ by omitting the dependence in $N$) the quantization parameter:
\begin{equation}
\hbar_N=\frac{(b-a)\beta}{N}.
\label{hbar}
\end{equation}
Finally, we define the quantization map as follows
\begin{equation}
T_N:C^\infty(\mathrm{\Sigma})\rightarrow M_N(\mathbb{C})\qquad T_N(f)=\sum_{n,m=1}^N f_{n-m}(z(n,m))E_{n,m}
\label{quantmap}
\end{equation}
where $E_{n,m}$ is the standard basis of $M_N(\mathbb{C})$.
\begin{mex}[The unit sphere]
Consider the local parametrization of the unit sphere given by
\begin{align*}
x(z,\varphi)=\sqrt{1-z^2}\cos(\varphi)\\
y(z,\varphi)=\sqrt{1-z^2}\sin(\varphi)
\end{align*}
equipped with its induced metric as a subspace of $\mathbb{R}^3$ and its Poisson bracket $\lbrace\cdot,\cdot\rbrace_{\mathbb{S}^2}$. Then, the quantization map on the unit sphere $T_N:C^\infty(\mathbb{S}^2)\rightarrow M_N(\mathbb{C})$ applied to the local coordinates functions $(x,y,z)$ gives:
\begin{equation}
T_N(x)=\frac{1}{2}\sqrt{1-z(n,n+1)^2}E_{n,n+1} + \frac{1}{2}\sqrt{1-z(n+1,n)^2}E_{n+1,n}
\end{equation}
and 
\begin{equation}
T_N(y)=\frac{1}{2i}\sqrt{1-z(n,n+1)^2}E_{n,n+1} - \frac{1}{2i}\sqrt{1-z(n+1,n)^2}E_{n+1,n}
\end{equation}
as well as 
\begin{equation}
T_N(z)=z(n,n)E_{n,n}.
\end{equation}
More generally, the example of the unit sphere may be extended to an ellipsoid with implicit representation:
\begin{equation}
\frac{x^2}{a^2}+\frac{y^2}{b^2}+\frac{z^2}{c^2}=1
\end{equation}
where $a,b$ and $c$ are the length of the semi-axes. The local coordinate system is given by 
\begin{equation}
x(z,\varphi)=a\sqrt{1-z^2}\cos(\varphi),\qquad y(z,\varphi)=b\sqrt{1-z^2}\sin(\varphi),
\end{equation}
and the $z$ coordinate. The quantization of the ellipsoid using $T_N$ is then defined by the matrix coordinates:
\begin{align}
&T_N(x)=\frac{a}{2}\sqrt{1-z(n,n+1)^2}E_{n,n+1} + \frac{a}{2}\sqrt{1-z(n+1,n)^2}E_{n+1,n},\\
&T_N(y)=\frac{b}{2i}\sqrt{1-z(n,n+1)^2}E_{n,n+1} - \frac{b}{2i}\sqrt{1-z(n+1,n)^2}E_{n+1,n},\\
&T_N(z)=cz(n,n)E_{n,n}.
\end{align}
\end{mex}
\noindent
We are now ready to prove that the quantization map \eqref{quantmap} satisfies the axioms \eqref{ax1}-\eqref{ax3} of a matrix quantization.
\begin{mth}
Consider $(\mathrm{\Sigma},\omega)$ a compact toric surface embedded in $\mathbb{R}^3$ equipped with the induced metric tensor $g$. The map defined: 
\begin{equation*}
T_N:C^\infty(\mathrm{\Sigma})\rightarrow M_N(\mathbb{C})\qquad T_N(f)=\sum_{n,m=1}^N f_{n-m}(z(n,m))E_{n,m}
\end{equation*}
is a quantization map.
\end{mth}
\begin{proof}
Using Parseval's inequality on $f\in L^2(S^1)$, we get that 
\begin{equation}
\|T_N(f)\|\leq \|f\|_{L^2}
\end{equation}
and thus $\|T_N\|<1$ uniformly, which proves \eqref{ax1}. Condition \eqref{ax2} follows from the following calculation
\begin{equation*}
(T_N(f)T_N(g))_{n,m}=\sum_{1<\ell<N}f_{\ell-n}(z(n,\ell))g_{m-\ell}(z(\ell,m))
\end{equation*}
and looking at the Taylor expansion of the Fourier coefficients close to a value $z(n,m)$, we get
\begin{equation*}
(T_N(f)T_N(g))_{n,m}=\sum_{1<\ell<N}f_{\ell-n}(z(n,m)+(\ell-m)\delta z)g_{m-\ell}(z(n,m)+(\ell-n)\delta z)
\end{equation*}
with
\begin{equation}
\delta z = \frac{(b-a)\beta}{N}, \qquad \beta>0.
\end{equation}
Since $\widehat{f}$ is compactly supported, $f_n=0$ for all $|n|\geq \delta$, then we have 
\begin{equation}
(\ell - n)\delta z \leq \frac{\delta}{N},
\end{equation}
for $1<\ell<N$. And thus, using \eqref{fouriercoeff}, we get 
\begin{equation*}
(T_N(f)T_N(g))_{n,m}=(T_N(fg))_{n,m} + \delta z\left( T_N\left( \partial_\varphi f\partial_z g\right)  - T_N\left( \partial_z f\partial_\varphi g\right) \right)_{n,m} +O(\delta z^2)
\end{equation*}
from which \eqref{ax2} follows:
\begin{equation}
T_N(f)T_N(g) = T_N(fg)+\frac{i\beta}{N}T_N\left( \left\lbrace f,g \right\rbrace \right) + O\left( \frac{1}{N^2}\right). 
\label{approx}
\end{equation}
Finally, Condition \eqref{ax3} can be obtained using condition \eqref{ax2} and computing the commutator:
\begin{equation}
\left[ T_N(f),T_N(g)\right] =\frac{2i}{N}T_N\left( \left\lbrace f,g \right\rbrace \right)+O\left( \frac{1}{N^2}\right) 
\end{equation}
which concludes the proof.
\end{proof}
\begin{mdef}[Smooth matrix quantization]
A quantization map $T_N$ is said to be a \textit{smooth} matrix regularization of $(\mathrm{\Sigma},\omega)$ if it holds that for any $f,h\in C^\infty(\mathrm{\Sigma})$, there exists $A_k(f,g)\in C^\infty(\mathrm{\Sigma})$ such that
\begin{equation}
\frac{1}{i\hbar}\left[ T_N(f),T_N(g)\right]=\sum_kc_{k,\alpha}(f,g)T_N(A_k(f,g)) 
\label{smooth}
\end{equation}
for some $c_{k,\alpha}(f,g)\in \mathbb{R}$.
\end{mdef}
\begin{mprop}
The quantization map \eqref{quantmap} defines a \textit{smooth} matrix regularization of $(\mathrm{\Sigma},\omega)$.
\label{remainder}
\end{mprop}
\begin{proof}
The statement follows by improving the approximation \eqref{approx} using Taylor remainder theorem:
\begin{equation}
T_N(f)T_N(g) = T_N(fg)+\frac{i\beta}{N}T_N\left( \left\lbrace f,g \right\rbrace \right) -\frac{\beta^2}{N^2}T_N\left(\partial^2_\varphi f\partial^2_z g|_{\xi_z}+ \partial^2_z f|_{\eta_z}\partial^2_\varphi g \right) 
\label{Taylor}
\end{equation}
from which \eqref{smooth} follows by computing the commutator.
\end{proof}
\noindent
We can now show how the square-root of the determinant of the metric tensor $\sqrt{|g|}$, which contains curvature information about the surface, can be approximated by matrix commutators.
\begin{mprop}
\label{propmetric}
Consider $(\mathrm{\Sigma},\omega)$ a compact toric surface embedded in $\mathbb{R}^3$ equipped with the induced metric tensor $g$. Let $(x_1,x_2,x_3)$ the local coordinate embedding system. Then, the image of $\sqrt{|g|}$ under the quantization map \eqref{quantmap} satisfies
\begin{equation*}
T_N\left( \sqrt{|g|}\right) =\gamma_N + O\left( \frac{1}{N^2}\right) 
\end{equation*}
where 
\begin{equation*}
\gamma_N=\sqrt{\sum_{i>j}-\left( \frac{[X^i_\alpha,X^j_\alpha]}{\hbar_\alpha}\right)^2 }
\end{equation*}
such that $T_N(x_i)=X_i$ for $i\in \left\lbrace 1,2,3\right\rbrace $.
\end{mprop}
\begin{proof}
Let us recall Equation \eqref{determetric}, where the metric determinant is expressed in local coordinates:
\begin{equation}
|g|=-\sum_{i>j}\left\lbrace x_i,x_j \right\rbrace^2.
\end{equation}
Then, if we apply the quantization map to both sides of the equation above and using linearity, we get 
\begin{equation}
T_N(|g|)=-\sum_{i>j}T_N(\left\lbrace x_i,x_j \right\rbrace^2)
\end{equation}
and using conditions \eqref{ax2}-\eqref{ax3} of the quantization, ti follows that
\begin{equation}
T_N(|g|)=-\sum_{i>j}\left(\frac{\left[ X_i,X_j\right] }{\hbar} \right) ^2 + O\left(\frac{1}{N^2} \right).
\end{equation}
Again, from condition \eqref{ax2}, we obtain
\begin{equation}
T_N(|g|)=T_N\left( \sqrt{|g|}\right)^2 + O(\frac{1}{N^2})
\end{equation}
which allows us to conclude the proof.
\end{proof}
\begin{mlem}
\label{leminv}
Let $f\in C^\infty(\mathrm{\Sigma})$ be a nowhere vanishing function. Then,
\begin{equation*}
T_N\left( 1/f\right) =T_N(f)^{-1} + O\left( \frac{1}{N^2}\right). 
\end{equation*}
\end{mlem}
\begin{proof}
The results follows from condition \eqref{ax2} applied to $1/f$ and $f$.
\end{proof}
\noindent
One can also address the question of convergence of the trace of the matrix $T_N(f)$ for some $f\in C^\infty(\mathrm{\Sigma})$. 
\begin{mprop}
Consider $(\mathrm{\Sigma},\omega)$ a compact surface embedded in $\mathbb{R}^3$ with an $S^1$-axial symmetry. Let $\omega_g$ be the volume form on $\mathrm{\Sigma}$. Then, the quantization map \eqref{quantmap} satisfies:
\begin{equation*}
\lim_{N\rightarrow \infty}2\pi\hbar_N \mathrm{Tr}\ T_N(f) =\int_{\mathrm{\Sigma}}f\omega.
\end{equation*}
\label{limtrace}
\end{mprop}
\begin{proof}
By construction of the quantization map, we have
\begin{equation}
\mathrm{Tr}\ T_N(f) =\sum_{1\leq n\leq N} f_0(z(n,n))
\end{equation}
where $f_0$ is the first coefficient in the Fourier series of $f$. This can be re-expressed as
\begin{equation}
\mathrm{Tr}\ T_N(f) =\frac{1}{2\pi}\sum_{1\leq n\leq N} \int_{0}^{2\pi} f(z(n,n),\theta)d\theta
\end{equation}
thus, multiplying by $2\pi \hbar_N$, we recognize a Riemann summation on the $z$-coordinate. Hence,
\begin{equation}
\lim_{N\rightarrow\infty}2\pi\hbar_N \mathrm{Tr}\ T_N(f) =\int_{\mathrm{\Sigma}}f\omega_g
\end{equation}
which is the desired property.
\end{proof}
\begin{mrmk}
A necessary condition for the limit \eqref{limtrace} to be satisfied is obtained by substituing $f\equiv 1$ in the equation. In this case, we see that $\beta$ must be:
\begin{equation}
\beta = \frac{\mathrm{vol}_g(\mathrm{\Sigma})}{2\pi(b-a)}
\end{equation}
for \eqref{limtrace} to hold.
\end{mrmk}
\begin{mth}
\label{th1}
Consider $(\mathrm{\Sigma},\omega)$ a compact toric surface embedded in $\mathbb{R}^3$ equipped with the induced metric tensor $g$. Let $(x_1,x_2,x_3)$ be the local coordinate embedding system. Denote by $\mathrm{\Delta_g}:C^\infty(\mathrm{\Sigma})\rightarrow C^\infty(\mathrm{\Sigma})$ the Laplace-Beltrami operator on $(\mathrm{\Sigma},\omega)$. Then, we have the following approximation of the Laplace-Beltrami operator through the quantization map:
\begin{equation}
 T_N(\mathrm{\Delta}_g f)=\mathrm{\Delta}_N T_N(f) + O\left( \frac{1}{N^2}\right) 
\end{equation}
where
\begin{equation}
\mathrm{\Delta}_N = -\frac{1}{\hbar^2}\gamma^{-1}_N [X_i,\gamma^{-1}_N[X_i,\cdot ] ]
\end{equation}
is the noncommutative Laplace operator and $T_N(x_i)=X_i$ for $i\in \left\lbrace 1,2,3\right\rbrace $.
\end{mth}
\begin{proof}
Similarly to the metric determinant function, we recall the expression of the Laplace-Beltrami operator in local coordinates:
\begin{equation}
\mathrm{\Delta}_g(f)=\sum_i\frac{1}{\sqrt{|g|}}\left\lbrace x^i,\frac{1}{\sqrt{|g|}}\left\lbrace x^i, f \right\rbrace \right\rbrace. 
\label{lap2}
\end{equation}
Moreover, let us recall that $T_N$ is a smooth regularization according to Proposition \ref{remainder} and we adapt Equation \eqref{approx}, we get
\begin{equation}
\frac{1}{i\hbar}\left[ T_N(x^i),T_N(f)\right]=T_N(\lbrace x^i,f\rbrace)-\frac{1}{\hbar^2}T_N(A(x^i,f)). 
\end{equation}
Then applying the quantization map $T_N$ to Equation \eqref{lap2} and using conditions \eqref{ax1}-\eqref{ax3}, we deduce that
\begin{equation}
T_N\mathrm{\Delta}_g(f)=\sum_i T_N(|g|^{-\frac{1}{2}})\left[ T_N(x^i),T_N(|g|^{-\frac{1}{2}})\left[ T_N(x^i),T_N(f)\right] \right] +O\left(\frac{1}{N^2}\right). 
\end{equation}
Finally, using Lemma \ref{leminv} and Proposition \ref{propmetric}, we get the desired result.
\end{proof}
\begin{mdef}[Noncommutative Laplacian]
Consider $(\mathrm{\Sigma},\omega)$ a compact toric surface embedded in $\mathbb{R}^3$ equipped with the induced metric tensor $g$. Let $(x_1,x_2,x_3)$ the local coordinate embedding system. Let $T_N$ be the quantization map \eqref{quantmap} on $(\mathrm{\Sigma},\omega)$. The \textit{noncommutative Laplace operator} on $\mathrm{\Sigma}$ is the operator 
\begin{equation}
\mathrm{\Delta}_N : M_N(\mathbb{C})\rightarrow M_N(\mathbb{C})\qquad \mathrm{\Delta}_N(F) = -\frac{1}{\hbar^2}\gamma^{-1}_N [X_i,\gamma^{-1}_N[X_i,F] ].
\end{equation}
An \textit{eigenmatrix sequence} of $\mathrm{\Delta}_N$ is a convergent sequence $\left\lbrace F_N \right\rbrace $ such that there exists $f\in C^\infty(\mathrm{\Sigma})$ satisfying $T_N(f)=F_N$ and
\begin{equation*}
\mathrm{\Delta}_N(F_N)=\lambda_NF_N, \ \text{for all $N$ with,} \ \lim_{N\rightarrow \infty}\lambda_N=\lambda.
\end{equation*}
\end{mdef}
\noindent
We are now ready to prove a convergence theorem of the spectrum of the noncommutative Laplacian to the spectrum of thee Laplace-Beltrami operator.
\begin{mth}
Let $\left\lbrace F_N\right\rbrace $ be an eigenmatrix sequence of $\mathrm{\Delta}_N$. Then, $\left\lbrace F_N\right\rbrace $ converges to an eigenfunction of the Laplace-Beltrami $\mathrm{\Delta}_g$ on $\mathrm{\Sigma}$ with eigenvalue given by the limit $\lambda=\lim_{N\rightarrow\infty}\lambda_N$.
\end{mth}
\begin{proof}
We need to prove that $\lim_{N\rightarrow \infty }\|T_N(\mathrm{\Delta}f -\lambda f)\|=0$. We have the following inequality:
\begin{equation}
\|T_N(\mathrm{\Delta}f-\lambda f)\|\leq \|T_N(\mathrm{\Delta}f)-\mathrm{\Delta}_NT_N(f)\| + \|\mathrm{\Delta}_NT_N(f)-\lambda T_N(f)\|
\end{equation}
and after taking the limit and using Theorem \ref{th1}, we get that
\begin{equation}
\lim_{N\rightarrow \infty}\|T_N(\mathrm{\Delta}f-\lambda f)\|=\lim_{N\rightarrow \infty} \|\mathrm{\Delta}_NT_N(f)-\lambda T_N(f)\|
\end{equation}
Moreover, we have that 
\begin{equation}
\|\mathrm{\Delta}_NT_N(f)-\lambda T_N(f)\| \leq \|\mathrm{\Delta}_NT_N(f)-\lambda_N T_N(f)\|+|\lambda-\lambda_N|\|T_N(f)\|
\end{equation}
from which the result follows after taking the limit.
\end{proof}
\section{Quantization of higher genus surfaces}
\label{Section4}
Up to this point, we have only considered a manifold with a unique coordinate system, i.e. with a global embedding coordinate system in $\mathbb{R}^3$. In this section, we treat the general case of a manifold given by several coordinate charts. In order to do so, first notice that for two disjoint manifolds (embedded in the same space), the algebra of functions on the union is the direct sum of the function spaces of the manifolds. Thus, in order to derive matrix regularizations in the general case, we first present results on basic operations on matrix systems and quantization of matrix valued functions. We then explain how the quantization of two sets of coordinate charts are glued together to produce a matrix quantization of the entire manifold.
\subsection{Operations on matrix regularization}
We start this section with the unitary transformation of matrix regularizations. Let $(\mathrm{\Sigma}_i,\omega_i)$ be two toric surfaces with embedding coordinates $(x_i,y_i,z_i)$ into $\mathbb{R}^3$ for $i=1,2$. The embedding coordinates of the disjoint union $\mathrm{\Sigma}_1\sqcup \mathrm{\Sigma}_2$ are given by
\begin{equation*}
q_1=\left( 
\begin{array}{cc}
x_1(\theta,z) & 0 \\ 
0 & x_2(\theta,z)
\end{array} 
\right), \quad 
q_2=\left( 
\begin{array}{cc}
y_1(\theta,z) & 0 \\ 
0 & y_2(\theta,z)
\end{array} 
\right), \quad
q_3=\left( 
\begin{array}{cc}
z& 0 \\ 
0 & z
\end{array} 
\right).
\end{equation*}
We first extend the quantization map $T_N$ defined in \eqref{quantmap} to matrix valued functions,
\begin{equation}
T_N:C(S^1,M_2(\mathbb{C}))\rightarrow M_N(\mathbb{C})
\end{equation}
as follows: fix a value of $z$ and let $f\in C(S^1,M_2(\mathbb{C}))$ such that 
\begin{equation}
f=\left( 
\begin{array}{cc}
f_{11}(\theta,z) & f_{12}(\theta,z) \\ 
f_{21}(\theta,z) & f_{22}(\theta,z)
\end{array} 
\right). 
\end{equation}
The Fourier expansions of the coefficients of $f$ are given by
\begin{equation}
f_{ij}(\theta)=\sum_{n\in\mathbb{Z}}\hat{f}^{ij}_n(z)e^{in\theta},\qquad 1\leq i,j\leq 2
\end{equation}
and we extend $T_N$ as follows
\begin{equation}
T_N(f)=\sum_{ij}T_N(f_{ij})\otimes E_{ij}\quad \text{and}\quad T_N(f_{ij})=\sum_{n,m}\hat{f}^{ij}_{n-m}(z(n,m))E_{n,m}
\label{extend}
\end{equation}
for $1\leq i,j\leq 2$.
\begin{mex}[Direct sum]
Let us assume that we have two matrix regularization $F_i$, which are regularizations of the functions $(\theta,z)\mapsto f_i(z,\theta)$, $i=1,2$. On one hand, the direct sum of the regularizations $F_1\oplus F_2$ denoted by $F$ is then given by:
\begin{equation}
F=\left( \begin{array}{cc}
F_1 & 0 \\ 
0 & F_2
\end{array} \right) =\sum_{1\leq n,m\leq N}\left(
\begin{array}{cc}
\hat{f}^{11}_{n-m}E_{n,m} & 0 \\ 
0 & \hat{f}^{22}_{n-m}E_{n,m}.
\end{array} 
 \right). 
\end{equation}
On the other hand, when applying the map $T_N$ to the function $f=f_1\oplus f_2$, we get:
\begin{equation}
T_N(f)=\sum_{1\leq n,m\leq N}\left( 
\begin{array}{cc}
\hat{f}^{11}_{n-m}& 0 \\ 
0 & \hat{f}^{22}_{n-m}
\end{array} 
\right) E_{n,m},
\end{equation}
therefore, for $P$ a suitable permutation operator, we notice that 
\begin{equation}
F=PT_N(f)P
\end{equation}
where the inner and outer indices of summations are exchanged. Thus, up to a unitary transformation, we can consider the operator $F$ and $T_N(f)$ to be equivalent.
\end{mex}
\noindent
Moreover, consider the matrix regularizations
\begin{equation}
T_N(x_i)=X_i,\quad T_N(y_i)=Y_i,\quad T_N(z_i)=Z_i
\end{equation}
of the set of coordinates $(x_i,y_i,z_i)$ of $(\mathrm{\Sigma}_i,\omega_i)$ for $i=1,2$. Then, we have the following result.
\begin{mlem}
Let $\lbrace T^i_{N}\rbrace$ be a matrix regularization of the compact toric manifold $(\mathrm{\Sigma}_i,\omega_i)$ for $i=1,2$. Then $\lbrace T^1_{N}\oplus T^2_N\rbrace$ is a matrix regularization of the disjoint union $(\mathrm{\Sigma}_1\sqcup \mathrm{\Sigma}_2,\omega_1\oplus \omega_2)$. 
\end{mlem}
\begin{proof}
The lemma follows straightforwardly from the fact that the algebra of smooth functions on $\mathrm{\Sigma}_1\sqcup \mathrm{\Sigma}_2$ can be identified with the direct sum $C^\infty(\mathrm{\Sigma})\oplus C^\infty(\mathrm{\Sigma})$.
\end{proof}
\noindent
In particular, the direct sum $\lbrace X_1\oplus X_2,Y_1\oplus Y_2, Z_1\oplus Z_2\rbrace$ is a matrix regularization of the coordinate functions $(x_1\oplus x_2,y_1\oplus y_2,z_1\oplus z_2)$ of $(\mathrm{\Sigma}_1\sqcup \mathrm{\Sigma}_2,\omega_1\oplus \omega_2)$. \\

\noindent
More generally, we can apply unitary transformation to local coordinates. Let $u$ be a $2\times 2$ unitary matrix with constant coefficients. Then, the unitary transformation of the coordinates $q_i$ is given by: 
\begin{equation}
\tilde{q}_i=u^*q_iu
\end{equation}
 and define the unitary operator
\begin{equation*}
U=\sum_{1\leq n\leq N} u E_{n,n}.
\end{equation*}
We notice that $U$ is equivalent the quantization of the function $u$ i.e. there exists a permutation matrix such that
\begin{equation}
P^*UP=T_N(u).
\end{equation}
In particular, using the extension \eqref{extend}, notice that $T_N(u)=U$ for some fixed $N$.
\begin{mcor}
Let $u$ be a unitary matrix and define the unitary operator
\begin{equation*}
T_N(u)=\sum_{1\leq i,j\leq 2} T_N(u_{ij})\otimes E_{ij}.
\end{equation*}
Then, the quantization of the unitary transformation $u^*fu$ of a periodic function $f$ is given by:
\begin{equation*}
T_N(u^*fu)=U^*T_N(f)U
\end{equation*}
where $U=P^*T_N(u)P$ for some permutation matrix $P$.
\end{mcor}
\noindent
In what follows, we identify a function $f$, respectively its quantization matrix $F$, with its unitary transformation $u^*fu$, respectively with $U^*FU$.
\begin{mex}[Direct sum of cylinders]
\label{directsum}
Consider the unitary transformation 
\begin{equation}
u=\frac{1}{\sqrt{2}}\left( 
\begin{array}{cc}
1 & -1 \\ 
1 & 1
\end{array} 
\right) 
\end{equation}
which mixes the diagonal and off-diagonal element of a given periodic matrix-valued function $f=\left(
\begin{array}{cc}
f_{1,1} & f_{1,2} \\ 
f_{2,1} & f_{2,2}
\end{array} 
 \right) $
 such that
 \begin{equation}
 u^*fu=\frac{1}{2}\left( 
 \begin{array}{cc}
 (f_{1,1}+f_{2,2})+(f_{1,2}+f_{2,1}) & (f_{2,2}-f_{1,1})+(f_{1,2}-f_{2,1}) \\ 
 (f_{2,2}-f_{1,1})-(f_{1,2}-f_{2,1}) & (f_{1,1}+f_{2,2})-(f_{1,2}+f_{2,1})
 \end{array} 
 \right) 
 \label{unittransfo}
 \end{equation}
 which are referred as "interlacing" surfaces in \cite{SHIMADA2004297, sykora_fuzzy_2017}. Such a transformation is fundamental when considering the quantization of the direct sum of two surfaces with (direct sum) coordinate functions having the following general form:
\begin{equation}
x=\left(
\begin{array}{cc}
x_1 & 0 \\ 
0 & -x_1
\end{array} 
\right), \quad 
y=\left( 
\begin{array}{cc}
y_1 & 0 \\ 
0 & -y_1
\end{array} 
\right), \quad
z=\left( 
\begin{array}{cc}
z_1 & 0 \\ 
0 & z_1
\end{array} 
\right). 
\end{equation}
In such case, the unitary transformation $u$ in Equation  \eqref{unittransfo} simplifies, and we get the equivalent coordinates
\begin{equation}
x'=\left(
\begin{array}{cc}
0 & x_1 \\ 
x_1 & 0
\end{array} 
\right), \quad 
y'=\left( 
\begin{array}{cc}
0 & y_1 \\ 
y_1 & 0
\end{array} 
\right), \quad
z'=\left( 
\begin{array}{cc}
z_1 & 0 \\ 
0 & z_1
\end{array} 
\right).
\end{equation}
The matrix regularization of $(x',y',z')$ is then given by 
\begin{equation}
T_N(x')=U^*T_N(x)U, \quad T_N(y')=U^*T_N(y)U, \quad T_N(z')=U^*T_N(z)U
\end{equation}
where the unitary operator $U$ is the quantization $T_N(u)$.
In particular, and as a concrete example in the case of the direct sum of two cylinder, one has:
\begin{equation}
x_1(\theta,z)=c_x(z)+r_x(z)\cos(\theta), \quad y_1(\theta,z)=c_y(z)+r_y(z)\cos(\theta)
\end{equation}
where $c_{x},c_y$ are center functions and $r_{x},r_y$ are radii.
\end{mex}
\subsection{Gluing of matrix regularizations}
\label{gluing}
In this section, we explain how the set of two quantizations of a toric manifold $\mathrm{\Sigma}$ defined by a set of two coordinate maps  $(x_i,y_i,z_i)$, $i=a,b$ can be glued together to form a matrix regularization of the whole manifold $\mathrm{\Sigma}$.\\
Let us assume that $(x_1,y_1,z)$ is are smooth local coordinates or $z\in [z_1,z_{n_2})$ and $(x_2,y_2,z)$ is are smooth local coordinates for $z\in (z_{n_2},z_N]$. Furthermore, if we let $r_i(\theta,z)=(x_i(\theta,z),y_i(\theta,z))$ for $i=a,b$, then let us assume in addition that we have the gluing condition:
\begin{equation}
\lim_{z\rightarrow z_{n_2}}r_a(\theta,z)=\lim_{z\rightarrow z_{n_2}}r_b(\theta,z) \ \text{and} \ \lim_{z\rightarrow z_{n_2}}\partial_z r_a(\theta,z)=\lim_{z\rightarrow z_{n_2}}\partial_z r_b(\theta,z).
\label{gluecond}
\end{equation}
The quantizations of the $z$ coordinate, on the interval $[z_1,z_{n2}]$, respectively on $[z_{n2},z_N]$, are given by the diagonal matrices:
\begin{equation*}
Z=\left( 
\begin{array}{ccccc}
z_1 &  &  &  &  \\ 
 & \ddots &  &  &  \\ 
 &  & z_{n_1} &  &  \\ 
 &  &  & \ddots &  \\ 
 &  &  &  & z_{n_2}
\end{array} 
\right) \qquad Z'=\left( 
\begin{array}{ccccc}
z_{n_2} &  &  &  &  \\ 
 & \ddots &  &  &  \\ 
 &  & z_{n_3} &  &  \\ 
 &  &  & \ddots &  \\ 
 &  &  &  & z_{N}
\end{array} 
\right).
\end{equation*}
Moreover, let us denote by $X_a$ the quantization of $x(\theta,z)$ for $z\in [z_1,z_{n_2}]$, such that is given by the band-matrix:
\begin{equation*}
X_a=\left( 
\begin{array}{cccccccc}
x^a_{1,1} & \cdots & x^a_{1,q} &  &  &  &  \\ 
\vdots & \ddots &  & \ddots & &  &  \\ 
x^a_{q,1} & & \ddots & &x^a_{n_1-q,n_1+q} & &  \\ 
& \ddots &  & x^a_{n_1,n_1}  & & \ddots &  \\ 
 &  & x^a_{n_1+q,n_1-q} &  & \ddots & & x^a_{n_2-q,n_2} \\ 
 &  & & \ddots & & \ddots & \vdots  \\ 
 &  &  & & x^a_{n_2,n_2-q} & \cdots & x^a_{n_2,n_2}
\end{array} 
\right)
\end{equation*}
where the coefficients are given by the Fourier series coefficient of the function $x_a$. Moreover, $X_b$ is the quantization of the function $x_b(\theta,z)$ for $z\in [z_{n_2},z_N]$, such that
\begin{equation*}
X_b=\left( 
\begin{array}{cccccccc}
x^b_{1,1} & \cdots & x^b_{1,q} &  &  &  &  \\ 
\vdots & \ddots &  & \ddots & &  &  \\ 
x^b_{q,1} & & \ddots & &x^b_{n_3-q,n_3+q} & &  \\ 
& \ddots &  & x^b_{n_3,n_3}  & & \ddots &  \\ 
 &  & x^b_{n_3+q,n_3-q} &  & \ddots & & x^b_{N-q,N} \\ 
 &  & & \ddots & & \ddots & \vdots  \\ 
 &  &  & & x^b_{N,N-q} & \cdots & x^b_{N,N}
\end{array} 
\right). 
\end{equation*}
Then, the matrix quantization $X$ of the direct sum $x_a\oplus x_b$ is obtained by concatenating the matrices $X_a$ and $X_b$ along the diagonals:
\begin{equation*}
X=\left( 
\begin{array}{cccccccccc}
x^a_{1,1} & \cdots & x^a_{1,q} &  &  &  & &\\ 
\vdots & \ddots &  & \ddots & &  & &\\ 
x^a_{q,1} & & \ddots & &x^a_{n_1-q,n_1+q} & & &\\ 
& \ddots &  & x^a_{n_1,n_1}  & & x^b_{n_3-q,n_3+q}  & &\\ 
 &  & x^a_{n_1+q,n_1-q} &  & x^b_{n_3,n_3} & & \ddots &\\ 
 &  &  & x^b_{n_3+q,n_3-q}  & & \ddots & & x^b_{N-q,N}\\ 
 &  &  & & \ddots & &\ddots & \vdots \\
 &  &  & & & x^b_{N,N-q}  & \cdots&x^b_{N,N} 
\end{array} 
\right).
\end{equation*}
The compatibility condition giving $X=T_N(x_a\oplus x_b)$ i.e. a matrix regularization of $x_a\oplus x_b$, follows from the gluing condition \eqref{gluecond} which insures that the Fourier coefficients of $X_a$ smoothly converge to the coefficient of $X_b$  as $z$ goes to $z_{n_2}$. Similarly, the quantization of $y$-coordinate is obtained by gluing $Y_a=T_N(y_a)$ and $Y_b=T_N(y_b)$ along the diagonals.
\section{Numerical results}
\label{Section5}
The numerical results we present in this section intend to demonstrate the effectiveness of the matrix Laplacian in approximating the spectrum of the Laplace-Beltrami operator on compact Riemann surfaces. The method relies on the MATLAB sparse eigenvalues solver \texttt{eigs} to compute the first set of  eigenvalues of the matrix Laplacian.
\subsection{The unit sphere}
 We start with the unit sphere, given in cartesian coordinates $x,y,z$:
 \begin{equation}
 x^2+y^2+z^2=1.
 \end{equation}
The analytic eigenvalues are explicitly known: the $k$-th eigenvalue is given by $\lambda_k=k(k+1)$ with the eigenspace $H_k$ having dimension
 \begin{equation}
\mathrm{dim}\ H_{\lambda_k}=\frac{(n+k)!}{n!k!}-\frac{(n+k-2)!}{n!(k-2)!}.
\end{equation}
We represent the sphere parametrically in cylindrical coordinates:
\begin{equation}
 x=\sqrt{1-z^2}\cos\theta,\quad y= \sqrt{1-z^2}\sin\theta,
\end{equation}
where $-1<z<1$ and $-\pi<\theta<\pi$.\\
 In Table \ref{tab2}, with see that the first eigenvalues of the noncommutative Laplace operator agree, with discrepancies on the order of $\hbar$ depending on the discretization parameter of the $z$-axis. This corroborates the theoretical convergence results of the matrix Laplace operator to the standard Laplace-Beltrami operator established in Theorem \eqref{th1}. It validates the noncommutative Laplace as a discrete approximation of the Laplace-Beltrami operator. The slight deviations observed on the eigenvalues are consistent with truncation effects inherent in finite-dimensional approximations, and they diminish as the size of the matrices increases.
\begin{table}[htp!]
\begin{center}
\begin{tabular}{|c|c|c|}
\hline 
$n$ & Analytic eigenvalues & Numerical eigenvalues $\hbar=0.001$\\ 
\hline 
$0$ & $0$ & $3.130828929442941\times 10^{-12}$ \\ 
\hline 
$1$ & $2$ & $2.000012087392884$ \\ 
\hline 
 & $2$ & $2.000023382754058$ \\ 
\hline 
& $2$ & $2.000089698069393$ \\ 
\hline 
$2$ & $6$ & $6.000039448443108$ \\ 
\hline 
 & $6$ & $6.000045822355212$ \\ 
\hline 
 & $6$ & $6.000074981785781$ \\ 
\hline 
 & $6$ & $6.000154377871050$ \\ 
\hline 
 & $6$ & $6.000400029350841$ \\ 
\hline 
\end{tabular} 
\end{center}
\caption{Eigenvalues of the Laplace-Beltrami operator on $\mathbb{S}^2$}
\label{tab2}
\end{table}
\subsection{The ellipsoid}
 In order to illustrate the flexibility of the method for geometries without high degrees of symmetry, we consider the more general case of an ellipsoid. Let $0<a<b$, the ellipsoid is given in cartesian coordinates $x,y,z$ by
\begin{equation}
\frac{x^2}{a^2}+ \frac{y^2}{a^2}+ \frac{z^2}{b^2}=1.
\end{equation}
We represent it in the parametric form by
\begin{equation}
 x=a\sqrt{1-z^2}\cos\theta,\quad y= a\sqrt{1-z^2}\sin\theta,
\end{equation}
where $-1<z<1$ and $-\pi<\theta<\pi$.\\
Unlike the sphere, the ellipsoid does not admit closed-form spectral data, but approximate eigenvalues can be obtained analytically \cite{volkmer_eigenvalues_2023, volkmer_laplace-beltrami_2024, eswarathasan_laplace-beltrami_2021}. In fact,
the study of Laplace’s equation on $\mathbb{R}^3$ using ellipsoidal-type coordinate systems is well-developed and is centered around the analysis of the Lamé equation \cite{DomínguezRiveraNigamOvall+2022+821+837}.\\ 
We consider the special  case of the ellipsoid for $a=\frac{1}{2}$ and $b=1$. Our computations of the eigenvalues are presented in Table \ref{tab1}. It shows an agreement of order $o(\hbar)$ with these references, indicating that the noncommutative Laplacian adapts naturally to deformations of the sphere. Importantly, the matrix Laplace operator remains a robust approximation under the anisotropy of the metric.
\begin{table}[htp!]
\begin{center}
\begin{tabular}{|c|c|}
\hline 
 Analytic eigenvalues & Numerical eigenvalues $\hbar=0.001$\\ 
\hline 
 $0$ & $8.554268404736831\times 10^{-13}$ \\ 
\hline 
 $9.4963551264$ & $9.497207579877351$ \\ 
\hline 
 $32.9870647190$  & $32.99092636681126$ \\ 
\hline 
 $70.0448683054$ & $70.02397573246742$ \\ 
\hline 
  $124.2351884130$ & $123.7890689580589$ \\ 
\hline 
  $244.4546288484$ & $245.2607627770743$ \\ 
\hline 
\end{tabular} 
\end{center}
\caption{Eigenvalues of the Lapalce-Beltrami operator on the ellipsoïd with $a=\frac{1}{2}$ and $b=1$.}
\label{tab1}
\end{table}
\noindent
These results warrant several important observations. First, the convergence rate of the scheme corroborates the theoretical results of previous sections. Second, the operator-theoretic framework ensures that symmetries, such as the rotational invariance, are respected at the discrete level. Third, the approach is well-suited to geometries defined parametrically or implicitly, since it does not rely on a mesh but only on the algebraic structure of the embedding.
\subsection{The standard torus}
\begin{figure}[htp!]
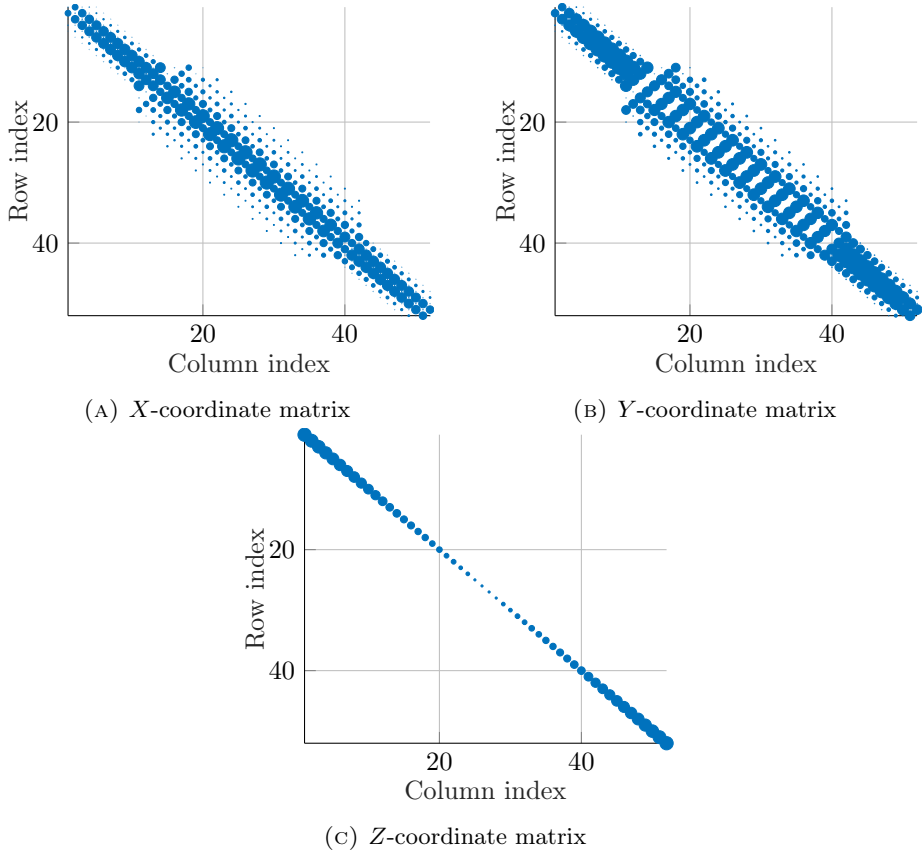

\centering
    \begin{subfigure}[t]{0.5\textwidth}
        \centering
        \tikzset{external/export next=true}
		\setlength\figureheight{0.200\textheight} 
		\setlength\figurewidth{0.8\linewidth} 
        \input{figures/CoordX.tikz} 
        \caption{$X$-coordinate matrix}
    \end{subfigure}%
    ~ 
    \begin{subfigure}[t]{0.5\textwidth}
        \centering
		\tikzset{external/export next=true}
		\setlength\figureheight{0.200\textheight} 
		\setlength\figurewidth{0.8\linewidth} 
        \input{figures/CoordY.tikz} 
        \caption{$Y$-coordinate matrix}
    \end{subfigure}
      ~ 
    \begin{subfigure}[t]{0.5\textwidth}
        \centering
		\tikzset{external/export next=true}
		\setlength\figureheight{0.200\textheight} 
		\setlength\figurewidth{0.8\linewidth} 
%
%
\definecolor{mycolor1}{rgb}{0.00000,0.44700,0.74100}%
\begin{tikzpicture}

\begin{axis}[%
width=0.951\figurewidth,
height=\figureheight,
at={(0\figurewidth,0\figureheight)},
scale only axis,
xmin=1,
xmax=52,
xlabel style={font=\color{white!15!black}},
xlabel={Column index},
y dir=reverse,
ymin=1,
ymax=52,
ylabel style={font=\color{white!15!black}},
ylabel={Row index},
axis background/.style={fill=white},
axis x line*=bottom,
axis y line*=left,
xmajorgrids,
ymajorgrids,
ylabel style={yshift=-5pt},xlabel style={yshift=2.5pt}
]
\addplot[scatter, only marks, mark=*, color=mycolor1, mark options={}, scatter/use mapped color={mark options={}, draw=mycolor1, fill=mycolor1}, visualization depends on={\thisrow{size} \as \perpointmarksize}, scatter/@pre marker code/.append style={/tikz/mark size=\perpointmarksize}] table[row sep=crcr]{%
x	y	size\\
1	1	2.5\\
2	2	2.43505383502854\\
3	3	2.36832733357835\\
4	4	2.2996655275195\\
5	5	2.22888957060432\\
6	6	2.15579124625644\\
7	7	2.08012573584461\\
8	8	2.00160192256359\\
9	9	1.91986911947076\\
10	10	1.83449846426336\\
11	11	1.74495610113028\\
12	12	1.65056322943381\\
13	13	1.44337567297406\\
14	14	1.44337567297406\\
15	15	1.32771251135949\\
16	16	1.32771251135949\\
17	17	1.20096115353815\\
18	18	1.20096115353815\\
19	19	1.0591481821704\\
20	20	1.0591481821704\\
21	21	0.895143592549291\\
22	22	0.895143592549291\\
23	23	0.693375245281536\\
24	24	0.693375245281536\\
25	25	0.400320384512718\\
26	26	0.400320384512718\\
27	27	0.400320384512718\\
28	28	0.400320384512718\\
29	29	0.693375245281536\\
30	30	0.693375245281536\\
31	31	0.895143592549291\\
32	32	0.895143592549291\\
33	33	1.0591481821704\\
34	34	1.0591481821704\\
35	35	1.20096115353815\\
36	36	1.20096115353815\\
37	37	1.32771251135949\\
38	38	1.32771251135949\\
39	39	1.44337567297406\\
40	40	1.44337567297406\\
41	41	1.65056322943381\\
42	42	1.74495610113028\\
43	43	1.83449846426336\\
44	44	1.91986911947076\\
45	45	2.00160192256359\\
46	46	2.08012573584461\\
47	47	2.15579124625644\\
48	48	2.22888957060432\\
49	49	2.2996655275195\\
50	50	2.36832733357835\\
51	51	2.43505383502854\\
52	52	2.5\\
};
\end{axis}
\end{tikzpicture}%
        \caption{$Z$-coordinate matrix}
    \end{subfigure}
    \caption{Quantized coordinates for the embedded torus $\mathbb{T}^2$}
 \label{fig4}   
\end{figure}

We consider a torus $\mathbb{T}^2$ with outer radius $R$ and inner radius, where 
\begin{equation*}
0<r<R.
\end{equation*} 
We represent the torus in $\mathbb{R}^3$ through its implicit equation
\begin{equation}
g(x,y,z)=(x^2+y^2+z^2+R^2-r^2)^2-4R^2(x^2+y^2)=0.
\label{eqTorus}
\end{equation}
An explicit local parametrization is also given by
\begin{equation}
x(\theta,z)=(R+r\sqrt{1-z^2})\cos\theta, \quad y(\theta,z)=(R+r\sqrt{1-z^2})\sin\theta
\end{equation}
using cylindrical coordinates $(\theta,z)$. A function $f:\mathbb{T}\rightarrow \mathbb{C}$ can be identified with a function $f:\mathbb{R}^2\rightarrow\mathbb{C}$ in cylindrical coordinates $(\theta,z)$. For a fixed value $z=z_0$, we get a restriction to a $2\pi$-periodic function $f:S^1\rightarrow\mathbb{C}$ on the curve obtained as a level set of \eqref{eqTorus}.\\
We consider the special the case of the standard torus with minor radius $r=1$ and major radius $R=2$. We find a foliation of the torus into level sets along the $z$-axis. A parametrization $r=(x(\theta),y(\theta))$ of a level set for a given value of $z$ is obtained by solving the gradient equation:
\begin{equation}
\frac{dr}{d\theta}=\mathrm{\nabla}g
\label{gradeq}
\end{equation}
for each value of $z$, using the implicit equation \eqref{eqTorus} and MATLAB \texttt{ode89} algorithm.\\
Consider the four critical point of the height Morse function that we denote by $(z_1,z_2,z_3,z_4)$. On the intervals $I_a=(z_1,z_2)$ and $I_c=(z_3,z_4)$, we obtain two U-shaped surface and on $I_b=(z_2,z_3)$ a direct sum of two-cylinders. Consider the local coordinates $(x_i,y_i,z_i)$, $i=a,b,c$ corresponding to each of these intervals. The matrix regularizations $(X,Y,Z)$ shown in Figure \ref{fig4} are obtained from the gluing of the three matrix regularization systems $(X_i,Y_i,Z_i)$, $i=a,b,c$, following Section \ref{gluing}. The matrices $(X_b,Y_b,Z_b)$ are the quantized coordinates of the direct sum of two cylinders and are obtained following Example \ref{directsum}. \\

\noindent
Our computations of the first non-zero eigenvalue of the matrix Laplacian are presented in Table \ref{tab3} for different values of $N$. It shows an agreement of order $o(\hbar)$ with analytical approximations of the first eigenvalue of the Laplace-Beltrami equation given in \cite{volkmer_laplace-beltrami_2021}. Importantly, these results show that the matrix Laplace operator remains a robust approximation for surfaces with higher genus.
\begin{table}
\begin{center}
\begin{tabular}{|c|c|c|}
\hline 
$N$ & Numerical eigenvalues & $|\lambda-\lambda_N|$ \\ 
\hline 
$17$ & $0.7722800134680878$ & $0.20445$ \\ 
\hline 
$47$ & $0.8837002729642458$ & $0.93031\times 10^{-1}$ \\ 
\hline 
$77$ & $0.9305686448517216$ & $0.46163\times 10^{-1}$ \\ 
\hline 
$107$ & $0.9546253109166054$ & $0.22106\times 10^{-1}$ \\ 
\hline 
$137$ & $0.9674492934177004$ & $0.92820\times 10^{-2}$ \\ 
\hline 
$167$ & $0.9737605063844723$ & $0.29708\times 10^{-2}$ \\ 
\hline 
$197$ & $0.9760115043650262$ & $0.71981\times 10^{-3}$ \\ 
\hline 
\end{tabular} 
\end{center}
\caption{First non-zero eigenvalue of the Laplace-Beltrami operator on the embedded torus $\mathbb{T}^2$ with $r=1$ and $R=2$.}
\label{tab3}
\end{table}
\subsection{The double torus}
Let us first recall thee general method to produce a polynomial equation whose level set is an $n$-torus in $\mathbb{R}^3$. Take
\begin{equation}
g(x,y,z)=(P(x)+y^2)^2+z^2-c^2
\end{equation}
where $c>0$, and $P$ is a polynomial of the form:
\begin{equation}
P(x)=a_{2k}x^{2k}+a_{2k-1}x^{2k-1}+\cdots+a_1x+a_0
\end{equation}
with $a_{2k}>0$. The surface $\mathrm{\Sigma}$ defined by the implicit equation:
\begin{equation}
g(x,y,z)=0
\end{equation}
is closed and bounded (since $P$ has even degree) hence compact. The surface $\mathrm{\Sigma}$ is a submanifold of $\mathbb{R}^3$ if and only if for each $p\in\mathrm{\Sigma}$, $d_pg\neq 0$  which is equivalent to requiring that the polynomials $P-c$ and $P+c$ have only simple roots. And, if the polynomial $P-c$ has exactly $2$ simple roots and the polynomial $P+c$ has exactly $2g$ simple roots, then $\chi(\mathrm{\Sigma})=2-2g$ and $\mathrm{\Sigma}$ is a surface of genus $g$.
\subsubsection*{Explicit construction of P} Let $g>0$
\begin{align}
&G(t)=(t-1)(t-2^2)\cdots(t-g^2)\quad \text{and}\quad M=\max_{0\leq t\leq g^2+1}G(t), \ \alpha\in (0,2c/M)\\
&Q(x)=\alpha G(x)-c\quad \text{and}\quad P(x)=Q(x^2).
\end{align}
One can see that $Q+c$ has exactly $g$ simple roots, hence $P+c$ has exactly $2g$ simple roots. For $t\in [0,g^2+1]$, the function $Q(t)-c$ has no zero. On the other hand, for $t\geq g^2+1$, $Q(t)-c$ is strictly growing and has exactly one zero. Consequently, the polynomial $P-c$ has exactly $2$ simple roots and the surface $\mathrm{\Sigma}$ defined above is a genus $g$ compact Riemann surface.\\

\noindent
We are considering the double torus defined by the implicit equation
\begin{equation}
g(x,y,z)=\left( \frac{1}{21}( z^2 - 1 )(z^2 - 4) + ( x^2 - 1 ) \right) ^2 +y^2 -1 =0
\end{equation}
The level sets of the double torus for several values of $z$ are shown in Figure \ref{slices} below. Similarly to the standard torus, these are obtained by solving Equation \eqref{gradeq} using the implicit equation \eqref{eqTorus} and MATLAB \texttt{ode89} algorithm.\
\begin{figure}[htp!]
\centering
		\tikzset{external/export next=true}
		\setlength\figureheight{0.200\textheight} 
		\setlength\figurewidth{0.8\linewidth} 
		\input{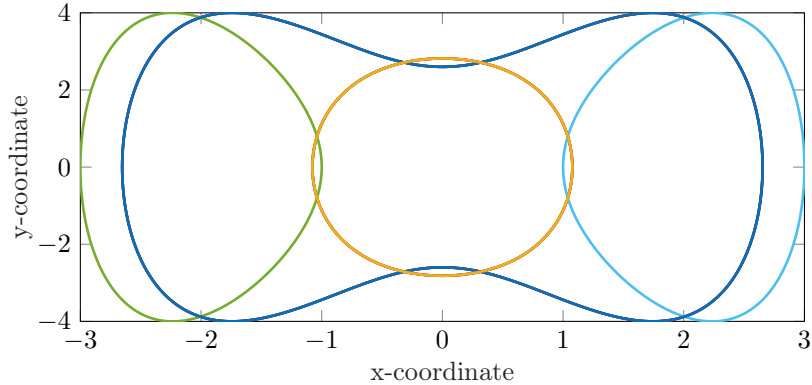}
\caption{Slices of the embedded torus for different values of $z$.}
\label{slices}
\end{figure}
Again, we consider the critical points of the height Morse function that we denote by $(z_1,z_2,z_3,z_4,z_5,z_6)$. On the intervals $I_a=(z_1,z_2)$, $I_e=(z_5,z_6)$, we obtain two U-shaped surface and on $I_b=(z_2,z_3)$ and $I_d=(z_4,z_5)$ a direct sum of two-cylinders. Consider the local coordinates $(x_i,y_i,z_i)$, $i=a,b,c$ corresponding to each of these intervals. The matrices $(X_b,Y_b,Z_b)$ and $(X_d,Y_d,Z_d)$ are the quantized coordinates of the direct sum of two cylinders and are obtained following Example \ref{directsum}. To the best of the authors’ knowledge, no closed-form expression or analytic approximation of the eigenvalues is available for the double torus.; nevertheless, we present Figure \ref{eigs2tor} that shows the first set of eigenvalues of the Laplace-Beltrami operator on the double torus for $\hbar=0.005$.
\begin{figure}[htp!]
\centering
\tikzset{external/export next=true}
\setlength\figureheight{0.200\textheight} 
\setlength\figurewidth{0.8\linewidth} 
%
%
\begin{tikzpicture}

\begin{axis}[%
width=0.951\figurewidth,
height=\figureheight,
at={(0\figurewidth,0\figureheight)},
scale only axis,
xmin=0,
xmax=30,
xlabel style={font=\color{white!15!black}},
xlabel={$k$ - wave number},
ymin=0,
ymax=30,
ylabel style={font=\color{white!15!black}},
ylabel={$\lambda_k$},
axis background/.style={fill=white},
ylabel style={yshift=-5pt},xlabel style={yshift=2.5pt}
]
\addplot [color=red, draw=none, mark=o, mark options={solid, red}, forget plot]
  table[row sep=crcr]{%
1	3.05866443284231e-13\\
2	0.781926783819231\\
3	3.25639197763099\\
4	6.2222172561926\\
5	8.3132849007176\\
6	9.05832804212302\\
7	9.368793401876\\
8	9.99339682821047\\
9	12.7193372389907\\
10	13.2129031155598\\
11	13.7864360463342\\
12	15.4383431302605\\
13	17.3413819242989\\
14	17.7894031781488\\
15	19.7371007366552\\
16	20.3829466703566\\
17	21.1647053862193\\
18	21.7753354130346\\
19	22.0757745367074\\
20	22.2237409530627\\
21	23.0455945875239\\
22	23.2207212366188\\
23	23.2416528163668\\
24	24.2554080077967\\
25	24.5110600775436\\
26	24.981058273186\\
27	26.0423296598021\\
28	26.6868186479616\\
29	26.8594961946403\\
30	27.1135440642116\\
};
\end{axis}
\end{tikzpicture}%
\caption{First set of eigenvalues of the Laplace-Beltrami operator of the embedded double torus.}
\label{eigs2tor}
\end{figure}
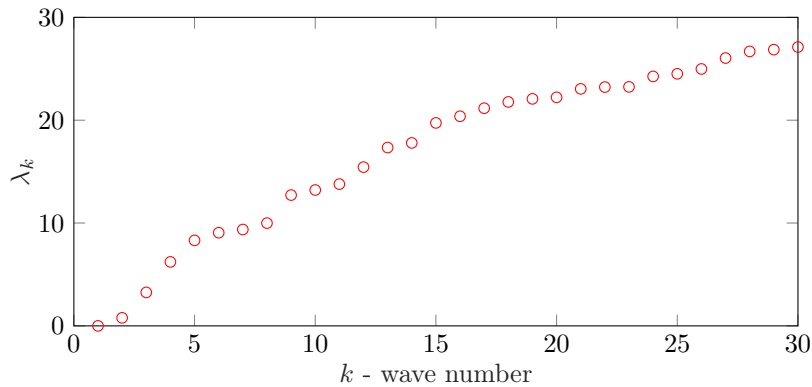
\section{Conclusion}
In this work, we have developed a noncommutative discretization of the Laplace Beltrami operator on compact Riemann surfaces. By replacing the function algebra of the surface with a sequence of finite-dimensional matrix algebras, we have introduced a noncommutative Laplacian defined via commutators of quantized coordinate functions. We have proved that this operator provides a consistent approximation to the classical Laplace–Beltrami operator, and we have demonstrated its effectiveness numerically on the sphere and surfaces of higher genus. The computed spectra illustrate both the convergence of eigenvalues and the robustness of the method under different geometries.\\

\noindent
This approach highlights the potential of structure-preserving discretizations that respect the underlying Poisson and symplectic structures of the surface. By embedding the problem into finite dimensional matrix spaces, we provide an alternative to mesh triangulations or point-based discretizations, that can capture geometric invariants in a finite-dimensional operator-theoretic setting.\\

\noindent
There remains several promising directions for further study. Extending the numerical tests to Riemann surfaces of higher genus would show the versatility of the method in more complex topological settings. A deeper numerical analysis of convergence rates and error estimates, including comparisons with finite element and spectral methods, would clarify its computational efficiency. It would also be natural to investigate the role of the noncommutative Laplacian in the study of time-dependent partial differential equations such as the heat or wave equation, thereby extending its applicability beyond static spectral problems. On the theoretical side, the relation between matrix regularization, Berezin–Toeplitz quantization, and fuzzy geometry suggests that the present approach could shed new light on the connection between quantization and spectral geometry.\\

\noindent
Beyond these mathematical questions, the method opens potential applications in several domains. In computational geometry and computer graphics, the intrinsic spectral information could be used for shape analysis and recognition without reliance on mesh quality. More broadly, the operator-theoretic discretization of the Laplacian may inspire new numerical schemes for partial differential equations on curved spaces, with implications ranging from fluid dynamics to general relativity. Since the construction does not fundamentally depend on triangulations, it may also be adapted to data-driven contexts where surfaces are represented only by point clouds. In mathematical physics, it provides a concrete framework for studying the spectra of quantized surfaces in membrane theory and string theory.\\

\nocite{*}
\bibliographystyle{amsplain}
\bibliography{References.bib}
\end{document}